\documentclass{article}
\usepackage{amssymb}
\usepackage{amsmath}
\usepackage{graphicx}

\newtheorem{theorem}{Theorem}

\newtheorem{definition}[theorem]{Definition}
\newtheorem{lemma}[theorem]{Lemma}

\newenvironment{proof}[1][Proof]{\noindent\textbf{#1.} }{\ \rule{0.5em}{0.5em}}

\begin{document}

\title{The $p$-adic Generalized Twisted $\left( h,q\right) $-Euler-$l$%
-Function and Its Applications}
\author{Mehmet Cenkci\\
Department of Mathematics, Akdeniz University, 07058-Antalya, Turkey \\
cenkci@akdeniz.edu.tr}
\date{}
\maketitle

\textbf{Abstract :} The main purpose of this paper is to construct
the $p$-adic twisted $\left(
h,q\right) $-Euler-$l$-function, which interpolates generalized twisted $%
\left( h,q\right) $-Euler numbers associated with a primitive
Dirichlet character $\chi $. This is a partial answer for the open
question which was remained in \cite{Kim-RimarXiv2006}. An
application of this function leads
general congruences systems for generalized twisted $\left( h,q\right) $%
-Euler numbers associated with $\chi $, in particular, Kummer-type
congruences for these numbers are obtained.



\section{\normalsize{Introduction}}

\hspace{0.15in}Let $\mathbb{N}$, $\mathbb{Z}$, $\mathbb{Q}$,
$\mathbb{R}$ and $\mathbb{C}$ denote, respectively, sets of
positive integer, integer, rational, real and complex numbers as
usual. Let $p$ be a fixed odd prime number and $x\in \mathbb{Q}$.
Then there exists integers $m$, $n$ and $\nu _{p}\left( x\right) $
such that $x=p^{\nu _{p}\left( x\right) }m/n$ and $p$ does not
divide either $m$ or $n$. Let $\left| \cdot \right| _{p}$ be
defined such that $\left| x\right| _{p}=p^{-\nu _{p}\left(
x\right) }$ and $\left|
0\right| _{p}=0$. Then $\left| \cdot \right| _{p}$ is a valuation on $%
\mathbb{Q}$ which satisfies the non-Archimedean property%
\begin{equation*}
\left| x+y\right| _{p}\leqslant \text{max}\left\{ \left| x\right|
_{p},\left| y\right| _{p}\right\} .
\end{equation*}%
\noindent Completion of $\mathbb{Q}$ with respect to $\left| \cdot \right|
_{p}$ is denoted by $\mathbb{Q}_{p}$ and called the field of $p$-adic
rational numbers. But $\mathbb{Q}_{p}$ itself is not complete with respect
to $\left| \cdot \right| _{p}$. $\mathbb{C}_{p}$ is the completion of the
algebraic closure of $\mathbb{Q}_{p}$ and $\mathbb{Z}_{p}=\left\{ x\in
\mathbb{Q}_{p}:\left| x\right| _{p}\leqslant 1\right\} $ is called the ring
of $p$-adic rational integers (see \cite{Koblitz1977}, \cite{Robert2000}).

Let $d$ be a fixed positive odd integer and let%
\begin{eqnarray*}
\mathbb{X} &=&\mathbb{X}_{d}=\underset{N}{\underleftarrow{\text{lim}}}\left(
\mathbb{Z}/dp^{N}\mathbb{Z}\right) \text{, }\mathbb{X}_{1}=\mathbb{Z}_{p}, \\
\mathbb{X}^{\ast } &=&\bigcup\limits_{\underset{\left( a,p\right) =1}{0<a<dp}%
}\left( a+dp^{N}\mathbb{Z}_{p}\right) , \\
a+dp^{N}\mathbb{Z}_{p} &=&\left\{ x\in \mathbb{X}:x\equiv a\left( \text{mod}%
dp^{N}\right) \right\} ,
\end{eqnarray*}%
\noindent where $N\in \mathbb{N}$ and $a\in \mathbb{Z}$ with $0\leqslant
a<dp^{N}$ (\cite{Cenkci-Can-Kurt2004}, \cite{Kim1999}, \cite%
{Kim-Jang-Rim-Pak2003}, \cite{Satoh1989}).

When talking about $q$-extensions, $q$ can variously be considered as an
indeterminate, a complex number $q\in \mathbb{C}$ or a $p$-adic number $q\in
\mathbb{C}_{p}$. If $q\in \mathbb{C}$, we normally assume that $\left|
q\right| <1$. If $q\in \mathbb{C}_{p}$, we assume that $\left| 1-q\right|
_{p}<p^{-1/\left( p-1\right) }$ so that for $\left| x\right| _{p}\leqslant 1$%
, we have $q^{x}=$exp$\left( x\text{log}q\right) $ (\cite%
{Cenkci-Can-Kurt2004}, \cite{Cenkci-Can2006}, \cite{Kim1999}, \cite{Kim2002a}%
, \cite{Kim2002b}). We use the notations%
\begin{equation*}
\left[ x\right] _{q}=\frac{1-q^{x}}{1-q}\text{ and }\left[ x\right] _{-q}=%
\frac{1-\left( -q\right) ^{x}}{1+q}\text{.}
\end{equation*}

We say that $f$ is uniformly differentiable function at a point $a\in
\mathbb{X}$, and denote this property by $f\in UD\left( \mathbb{X}\right) $,
if the quotient of the differences%
\begin{equation*}
F_{f}=\frac{f\left( x\right) -f\left( y\right) }{x-y}
\end{equation*}%
\noindent has a limit $l=f^{\prime }\left( a\right) $ as $\left( x,y\right)
\rightarrow \left( a,a\right) $. For $f\in UD\left( \mathbb{X}\right) $, the
$p$-adic invariant $q$-integral on $\mathbb{X}$ was defined by%
\begin{equation*}
I_{q}\left( f\right) =\int\limits_{\mathbb{X}}f\left( t\right) d\mu
_{q}\left( t\right) =\underset{N\rightarrow \infty }{\text{lim}}\frac{1}{%
\left[ dp^{N}\right] _{q}}\sum_{a=0}^{dp^{N}-1}f\left( a\right) q^{a}
\end{equation*}%
\noindent (cf. \cite{Kim1999}, \cite{Kim2002b}), where for any positive
integer $N$%
\begin{equation*}
\mu _{q}\left( a+dp^{N}\mathbb{Z}_{p}\right) =\frac{q^{a}}{\left[ dp^{N}%
\right] _{q}}
\end{equation*}%
\noindent (cf. \cite{Kim1999}, \cite{Kim2002a}, \cite{Kim2002b}).

The concept of \textit{twisted} has been applied by many authors to certain
functions which interpolate certain number sequences. In \cite{Koblitz1979},
Koblitz defined twisted Dirichlet $L$-function which interpolates twisted
Bernoulli numbers in the field of complex numbers. In \cite{Simsek2003},
Simsek constructed a $q$-analogue of the twisted $L$-function interpolating $%
q$-twisted Bernoulli numbers. Kim et.al. \cite{Kim-Jang-Rim-Pak2003} derived
a $p$-adic analogue of the twisted $L$-function by using $p$-adic invariant
integrals. By using the definition of $h$-extension of $p$-adic $q$-$L$%
-function which is constructed by Kim \cite{KimarXivb}, Simsek \cite%
{Simsek2006a,Simsek2006b} and Jang \cite{Jang2007} defined twisted $p $-adic
generalized $\left( h,q\right) $-$L$-function. In \cite{Satoh1989}, Satoh
derived $p$-adic interpolation function for $q$-Frobenius-Euler numbers.
Simsek \cite{Simsek2005} gave twisted extensions of $q$-Frobenius-Euler
numbers and their interpolating function $q$-twisted $l$-series. In \cite%
{Cenkci-Can-Kurt2004}, Cenkci et.al. constructed generalized $p$-adic
twisted $l$-function in $p$-adic number field. Recently, Kim and Rim \cite%
{Kim-RimarXiv2006} defined twisted $q$-Euler numbers by using $p$-adic
invariant integral on $\mathbb{Z}_{p}$ in the fermionic sense. In that
paper, they raised the following question: \textit{Find a }$p$\textit{-adic
analogue of the }$q$\textit{-twisted }$l$-\textit{function which
interpolates }$E_{n,\xi ,q,\chi }^{\left( h,1\right) }$\textit{, the
generalized twisted }$q$\textit{-Euler numbers attached to }$\chi $\textit{\ %
\cite{Kim2006}, \cite{Kim2007a}}. In a forthcoming paper, Rim et.al. \cite%
{Rim-Simsek-Kurt-Kim} answered this question by constructing partial $\left(
h,q\right) $-zeta function motivating from a method of Washington \cite%
{Washington1976,Washington1997}.

In this paper, we construct $p$-adic generalized twisted $\left( h,q\right) $%
-Euler-$l$-function by employing $p$-adic invariant measure on $p$-adic
number field. This is the answer of the part of the question posed in \cite%
{Kim-RimarXiv2006}. This way of derivation of $p$-adic generalized twisted $%
\left( h,q\right) $-Euler-$l$-function is different from that of \cite%
{Rim-Simsek-Kurt-Kim}, and leads an explicit integral representation for
this function. As an application of the derived integral representation, we
obtain general congruences systems for generalized twisted $q$-Euler numbers
associated with $\chi $, containing Kummer-type congruences.

\section{\normalsize{Generalized Twisted $q$-Euler Numbers}}

\hspace{0.15in}In this section, we give a brief summary of the concepts $p$-adic $q$%
-integrals and Euler numbers and polynomials. Let $UD\left( \mathbb{X}%
\right) $ be the set of all uniformly differentiable functions on $\mathbb{X}
$. For any $f\in UD\left( \mathbb{X}\right) $, Kim defined a $q$-analogue of
an integral with respect to a $p$-adic invariant measure in \cite%
{Kim1999,Kim2002b} which was called $p$-adic $q$-integral. The $p$-adic $q$%
-integral was defined as follows:%
\begin{equation*}
I_{q}\left( f\right) =\int\limits_{\mathbb{X}}f\left( t\right) d\mu
_{q}\left( t\right) =\underset{N\rightarrow \infty }{\text{lim}}\frac{1}{%
\left[ dp^{N}\right] _{q}}\sum_{a=0}^{dp^{N}-1}f\left( a\right) q^{a}.
\end{equation*}%
\noindent Note that%
\begin{equation*}
I_{1}\left( f\right) =\underset{q\rightarrow 1}{\text{lim}}I_{q}\left(
f\right) =\int\limits_{\mathbb{X}}f\left( t\right) d\mu _{1}\left( t\right) =%
\underset{N\rightarrow \infty }{\text{lim}}\frac{1}{dp^{N}}%
\sum_{a=0}^{dp^{N}-1}f\left( a\right)
\end{equation*}%
\noindent is the Volkenborn integral (see \cite{Robert2000}).

The Euler zeta function $\zeta _{E}\left( s\right) $ is defined by means of%
\begin{equation*}
\zeta _{E}\left( s\right) =2\sum_{k=1}^{\infty }\frac{\left( -1\right) ^{k}}{%
k^{s}}
\end{equation*}%
\noindent for $s\in \mathbb{C}$ with Re$\left( s\right) >1$ (cf. \cite%
{Kim2006}). For a Dirichlet character $\chi $ with conductor $d$, $d\in
\mathbb{N}$, $d$ is odd, the $l$-function associated with $\chi $ is defined
as (\cite{Kim2006})%
\begin{equation*}
l\left( s,\chi \right) =2\sum_{k=1}^{\infty }\frac{\chi \left( k\right)
\left( -1\right) ^{k}}{k^{s}}
\end{equation*}%
\noindent for $s\in \mathbb{C}$ with Re$\left( s\right) >1$. This function
can be analytically continued to whole complex plane, except $s=1$ when $%
\chi =1$; and when $\chi =1$, it reduces to Euler zeta function $\zeta
_{E}\left( s\right) $. In \cite{Kim2007b}, $\left( h,q\right) $-extension of
Euler zeta function is defined by%
\begin{equation*}
\zeta _{E,q}^{\left( h\right) }\left( s,x\right) =\left[ 2\right]
_{q}\sum_{k=0}^{\infty }\frac{\left( -1\right) ^{k}q^{hk}}{\left[ k+x\right]
_{q}^{s}}
\end{equation*}%
\noindent with $s$, $h\in \mathbb{C}$, Re$\left( s\right) >1$ and $x\neq $%
negative integer or zero. $\left( h,q\right) $-Euler polynomials are defined
by the $p$-adic $q$-integral as%
\begin{equation*}
E_{n,q}^{\left( h\right) }\left( x\right) =\int\limits_{\mathbb{X}}q^{\left(
h-1\right) t}\left[ t+x\right] _{q}^{n}d\mu _{-q}\left( t\right) ,
\end{equation*}%
\noindent for $h\in \mathbb{Z}$. \ $E_{n,q}^{\left( h\right) }\left(
0\right) =E_{n,q}^{\left( h\right) }$ are called $\left( h,q\right) $-Euler
numbers. In \cite{Kim2007b}, it has been shown that for $n\in \mathbb{Z}$, $%
n\geqslant 0$%
\begin{equation*}
\zeta _{E,q}^{\left( h\right) }\left( -n,x\right) =E_{n,q}^{\left( h\right)
}\left( x\right) ,
\end{equation*}%
\noindent thus we have%
\begin{equation*}
E_{n,q}^{\left( h\right) }\left( x\right) =\left[ 2\right]
_{q}\sum_{k=0}^{\infty }\left( -1\right) ^{k}q^{hk}\left[ k+x\right]
_{q}^{n},
\end{equation*}%
\noindent from which the following entails:%
\begin{equation*}
E_{n,q}^{\left( h\right) }\left( x\right) =\frac{\left[ 2\right] _{q}}{%
\left( 1-q\right) ^{n}}\sum_{j=0}^{n}\binom{n}{j}\left( -1\right) ^{j}q^{xj}%
\frac{1}{1+q^{h+j}}.
\end{equation*}

In \cite{Kim2006,Kim2007b}, $\left( h,q\right) $-extension of the $l$%
-function associated with $\chi $ is defined by%
\begin{equation*}
l_{q}^{\left( h\right) }\left( s,\chi \right) =\left[ 2\right]
_{q}\sum_{k=1}^{\infty }\frac{\chi \left( k\right) \left( -1\right)
^{k}q^{hk}}{\left[ k\right] _{q}^{s}}
\end{equation*}%
\noindent for $h$, $s\in \mathbb{C}$ with Re$\left( s\right) >1$. The
negative integer values of $s$ are determined explicitly by%
\begin{equation*}
l_{q}^{\left( h\right) }\left( -n,\chi \right) =E_{n,q,\chi }^{\left(
h\right) },
\end{equation*}%
\noindent for $n\in \mathbb{Z}$, $n\geqslant 0$ where $E_{n,q,\chi }^{\left(
h\right) }$ are the generalized $\left( h,q\right) $-Euler numbers
associated with $\chi $ defined by%
\begin{equation*}
E_{n,q,\chi }^{\left( h\right) } =\int\limits_{\mathbb{X}}\chi
\left( t\right) q^{\left( h-1\right) t}\left[ t\right]
_{q}^{n}d\mu _{-q}\left( t\right)\left( =\left[ 2\right]
_{q}\sum_{k=1}^{\infty }\chi \left( k\right) \left( -1\right)
^{k}q^{hk}\left[ k\right] _{q}^{n}\right) .
\end{equation*}

Now assume that $q\in \mathbb{C}_{p}$ with $\left| 1-q\right| _{p}<1$. From
the definition of $p$-adic invariant $q$ integral on $\mathbb{X}$, Kim \cite%
{Kim2006} defined the integral%
\begin{equation}
I_{-1}\left( f\right) =\underset{q\rightarrow -1}{\text{lim}}I_{q}\left(
f\right) =\int\limits_{\mathbb{X}}f\left( t\right) d\mu _{-1}\left( t\right)
\label{2.1}
\end{equation}%
\noindent for $f\in UD\left( \mathbb{X}\right) $. Note that%
\begin{equation}
I_{-1}\left( f_{1}\right) +I_{-1}\left( f\right) =2f\left( 0\right) ,
\label{2.2}
\end{equation}%
\noindent where $f_{1}\left( t\right) =f\left( t+1\right) $. Repeated
application of last formula yields%
\begin{equation}
I_{-1}\left( f_{n}\right) =\left( -1\right) ^{n}I_{-1}\left( f\right)
+2\sum_{j=0}^{n-1}\left( -1\right) ^{n-1-j}f\left( j\right) ,  \label{2.3}
\end{equation}%
\noindent with $f_{n}\left( t\right) =f\left( t+n\right) $.

Let $T_{p}=\bigcup\limits_{n\geqslant 1}C_{p^{n}}=\underset{n\rightarrow
\infty }{\text{lim}}\mathbb{Z}/p^{n}\mathbb{Z}$, where $C_{p^{n}}=\left\{
w\in \mathbb{X}:w^{p^{n}}=1\right\} $ is the cyclic group of order $p^{n}$.
For $w\in T_{p}$, let $\phi _{w}:\mathbb{Z}_{p}\rightarrow \mathbb{C}_{p}$
denote the locally constant function defined by $t\rightarrow w^{t}$.

For $f\left( t\right) =\phi _{w}\left( t\right) e^{zt}$, we obtain%
\begin{equation*}
\int\limits_{\mathbb{X}}\phi _{w}\left( t\right) e^{zt}d\mu _{-1}\left(
z\right) =\frac{2}{we^{z}+1}
\end{equation*}%
\noindent using (\ref{2.1}) and (\ref{2.2}), and%
\begin{equation*}
\int\limits_{\mathbb{X}}\chi \left( t\right) \phi _{w}\left( t\right)
e^{zt}d\mu _{-1}\left( t\right) =2\sum_{i=1}^{d}\frac{\chi \left( i\right)
\phi _{w}\left( i\right) e^{iz}}{w^{d}e^{dz}+1}
\end{equation*}%
\noindent using (\ref{2.1}) and (\ref{2.3}) (cf. \cite{Kim2006}). As a
consequence, the twisted Euler numbers and generalized twisted Euler numbers
associated with $\chi $ can respectively be defined by%
\begin{equation*}
\frac{2}{we^{z}+1}=\sum_{n=0}^{\infty }E_{n,w}\frac{z^{n}}{n!}\text{, and }%
2\sum_{i=1}^{d}\frac{\chi \left( i\right) \phi _{w}\left( i\right) e^{iz}}{%
w^{d}e^{dz}+1}=\sum_{n=0}^{\infty }E_{n,w,\chi }\frac{z^{n}}{n!},
\end{equation*}%
\noindent from which%
\begin{equation*}
\int\limits_{\mathbb{X}}t^{n}\phi _{w}\left( t\right) d\mu _{-1}\left(
t\right) =E_{n,w}\text{, and }\int\limits_{\mathbb{X}}\chi \left( t\right)
t^{n}\phi _{w}\left( t\right) d\mu _{-1}\left( t\right) =E_{n,w,\chi }
\end{equation*}%
\noindent follow.

Twisted extension of $\left( h,q\right) $-Euler zeta function is defined by%
\begin{equation*}
\zeta _{E,q,w}^{\left( h\right) }\left( s,x\right) =\left[ 2\right]
_{q}\sum_{k=0}^{\infty }\frac{\left( -1\right) ^{k}w^{k}q^{hk}}{\left[ k+x%
\right] _{q}^{s}}
\end{equation*}%
\noindent with $h$, $s\in \mathbb{C}$, Re$\left( s\right) >1$ and $x\neq $%
negative integer or zero. For $n\in \mathbb{Z}$, $n\geqslant 0$ and $h\in
\mathbb{Z}$, this function gives%
\begin{equation*}
\zeta _{E,q,w}^{\left( h\right) }\left( -n,x\right) =E_{n,q,w}^{\left(
h\right) }\left( x\right) ,
\end{equation*}%
\noindent where $E_{n,q,w}^{\left( h\right) }\left( x\right) $ are the
twisted $q$-Euler polynomials defined as%
\begin{equation*}
E_{n,q,w}^{\left( h\right) }\left( x\right)=\int\limits_{\mathbb{X}%
}q^{\left( h-1\right) t}\phi _{w}\left( t\right) \left[ x+t\right]
_{q}^{n}d\mu _{-q}\left( t\right)\left( =\left[ 2\right]
_{q}\sum_{k=0}^{\infty }\left( -1\right) ^{k}w^{k}q^{hk}\left[
k+x\right] _{q}^{n}\right)
\end{equation*}%
\noindent by using $p$-adic invariant $q$-integral on $\mathbb{X}$ in the
fermionic sense (cf. \cite{Kim-RimarXiv2006}, \cite{Rim-Simsek-Kurt-Kim}).
The following expressions for twisted $\left( h,q\right) $-Euler polynomials
can be verified from the defining equalities:%
\begin{eqnarray}
E_{n,q,w}^{\left( h\right) }\left( x\right) &=&\frac{\left[ 2\right] _{q}}{%
\left( 1-q\right) ^{n}}\sum_{j=0}^{n}\binom{n}{j}\left( -1\right) ^{j}q^{xj}%
\frac{1}{1+wq^{h+j}},  \label{2.4} \\
E_{n,q,w}^{\left( h\right) }\left( x\right) &=&\frac{\left[ 2\right] _{q}}{%
\left[ 2\right] _{q^{d}}}\left[ d\right] _{q}^{n}\sum_{a=0}^{d-1}q^{ha}w^{a}%
\left( -1\right) ^{a}E_{n,q^{d},w^{d}}^{\left( h\right) }\left( \frac{x+a}{d}%
\right) ,  \label{2.7}
\end{eqnarray}%
\noindent where $n$, $d\in \mathbb{N}$ with $d$ is odd. From (\ref{2.4}),
the twisted $\left( h,q\right) $-Euler polynomials can be determined
explicitly. A few of them are%
\begin{eqnarray*}
E_{0,q,w}^{\left( h\right) }\left( x\right) &=&\frac{1+q}{1+wq^{h}}, \\
E_{1,q,w}^{\left( h\right) }\left( x\right) &=&\frac{1+q}{1-q}\left( \frac{1%
}{1+wq^{h}}-\frac{q^{x}}{1+wq^{h+1}}\right) , \\
E_{2,q,w}^{\left( h\right) }\left( x\right) &=&\frac{1+q}{\left( 1-q\right)
^{2}}\left( \frac{1}{1+wq^{h}}-\frac{2q^{x}}{1+wq^{h+1}}+\frac{q^{2x}}{%
1+wq^{h+2}}\right) .
\end{eqnarray*}%
\noindent For $x=0$, $E_{n,q,w}^{\left( h\right) }\left( 0\right)
=E_{n,q,w}^{\left( h\right) }$ are called twisted $\left( h,q\right) $-Euler
numbers. Thus we can write%
\begin{equation*}
E_{n,q,w}^{\left( h\right) }\left( x\right) =\sum_{j=0}^{n}\binom{n}{j}q^{xj}%
\left[ x\right] _{q}^{n-j}E_{j,q,w}^{\left( h\right) }.
\end{equation*}

Let $\chi $ be a Dirichlet character of conductor $d$ with $d\in \mathbb{N}$
and $d$ is odd. Then the generalized twisted $\left( h,q\right) $-Euler
numbers associated with $\chi $ are defined as%
\begin{equation*}
E_{n,q,w,\chi }^{\left( h\right) }=\int\limits_{\mathbb{X}}\chi \left(
t\right) q^{\left( h-1\right) t}\phi _{w}\left( t\right) \left[ t\right]
_{q}^{n}d\mu _{-q}\left( t\right) .
\end{equation*}%
\noindent These numbers arise at the negative integer values of the twisted $%
\left( h,q\right) $-Euler-$l$-function which is defined by%
\begin{equation*}
l_{q,w}^{\left( h\right) }\left( s,\chi \right) =\left[ 2\right]
_{q}\sum_{k=1}^{\infty }\frac{\chi \left( k\right) \left( -1\right)
^{k}w^{k}q^{hk}}{\left[ k\right] _{q}^{s}}
\end{equation*}%
\noindent with $h$, $s\in \mathbb{C}$, Re$\left( s\right) >1$. Indeed, for $%
n\in \mathbb{Z}$, $n\geqslant 0$ and $h\in \mathbb{Z}$, we have%
\begin{equation*}
l_{q,w}^{\left( h\right) }\left( -n,\chi \right) =E_{n,q,w,\chi }^{\left(
h\right) }
\end{equation*}%
\noindent (cf. \cite{Kim-RimarXiv2006}, \cite{Rim-Simsek-Kurt-Kim}).

We conclude this section by stating the distribution property for
generalized twisted $\left( h,q\right) $-Euler numbers associated with $\chi
$, which will take a major role in constructing a measure in the next
section.

For $n$, $d\in \mathbb{N}$ with $d$ is odd, we have%
\begin{equation*}
E_{n,q,w,\chi }^{\left( h\right) }=\frac{\left[ 2\right] _{q}}{\left[ 2%
\right] _{q^{d}}}\left[ d\right] _{q}^{n}\sum_{a=1}^{d}q^{ha}w^{a}\chi
\left( a\right) \left( -1\right) ^{a}E_{n,q^{d},w^{d}}^{\left( h\right)
}\left( \frac{a}{d}\right) .
\end{equation*}

\section{\normalsize{$p$-adic Twisted $\left( h,q\right)
$-$l$-Functions}}

\hspace{0.15in}In this section we first focus on defining a
$p$-adic invariant measure, which is apparently an important tool
to construct $p$-adic twisted $\left(
h,q\right) $-Euler-$l$-function in the sense of $p$-adic invariant $q$%
-integral. We afterwards give the definition of $p$-adic twisted $\left(
h,q\right) $-Euler-$l$-function, together with Witt's type formulas for
twisted and generalized twisted $\left( h,q\right) $-Euler numbers.

Throughout, we assume that $\xi $ is the $r$th root of unity with $\left(
r,pd\right) =1$, where $p$ is an odd prime and $d$ is an odd natural number.
If $\left( r,pd\right) =1$, it has been known that $\left| 1-\xi \right|
_{p}\geqslant 1$ (see \cite{Koblitz1979}, \cite{Shiratani1985}) and $\xi $
lies in the cyclic group $C_{p^{n}}=\left\{ w:w^{p^{n}}=1\right\} $. The
following theorem plays a crucial role in constructing $p$-adic generalized
twisted $\left( h,q\right) $-Euler-$l$-function on $\mathbb{X}$.

\begin{theorem}
\label{thm1}Let $q\in \mathbb{C}_{p}$ with $\left| 1-q\right|
_{p}<p^{-1/\left( p-1\right) }$ and $\xi $ is the $r$th root of unity with $%
\left| 1-\xi \right| _{p}\geqslant 1$. For $N\in \mathbb{Z}$, $n\in \mathbb{Z%
}$, $n\geqslant 0$, let $\mu _{n,\xi ,q}^{\left( h\right) }$ be defined as%
\begin{equation*}
\mu _{n,\xi ,q}^{\left( h\right) }\left( a+dp^{N}\mathbb{Z}_{p}\right) =%
\left[ dp^{N}\right] _{q}^{n}\frac{\left[ 2\right] _{q}}{\left[ 2\right]
_{q^{dp^{N}}}}\left( -1\right) ^{a}\xi ^{a}q^{ha}E_{n,q^{dp^{N}},\xi
^{dp^{N}}}^{\left( h\right) }\left( \frac{a}{dp^{N}}\right) .
\end{equation*}%
\noindent Then $\mu _{n,\xi ,q}^{\left( h\right) }$ extends uniquely to a
measure on $\mathbb{X}$.
\end{theorem}

\begin{proof}
In order to show that $\mu _{n,\xi ,q}^{\left( h\right) }$ is a measure on $%
\mathbb{X}$, we need to show that it is a distribution and is bounded on $%
\mathbb{X}$.

To show it is a distribution on $\mathbb{X}$, we check the equality%
\begin{equation*}
\sum_{i=0}^{p-1}\mu _{n,\xi ,q}^{\left( h\right) }\left( a+idp^{N}+dp^{N+1}%
\mathbb{Z}_{p}\right) =\mu _{n,\xi ,q}^{\left( h\right) }\left( a+dp^{N}%
\mathbb{Z}_{p}\right) .
\end{equation*}%
\noindent Beginning the calculation from right hand side yields%
\begin{eqnarray*}
&&\hspace{-0.6in}\sum_{i=0}^{p-1}\mu _{n,\xi ,q}^{\left( h\right)
}\left(
a+idp^{N}+dp^{N+1}\mathbb{Z}_{p}\right)\\
&&\hspace{-0.5in}=\sum_{i=0}^{p-1}\left[ dp^{N+1}\right]
_{q}^{n}\frac{\left[ 2\right] _{q}}{\left[ 2\right]
_{q^{dp^{N+1}}}}\left( -1\right)
^{a+idp^{N}}\xi ^{a+idp^{N}}q^{h\left( a+idp^{N}\right) }E_{n,q^{dp^{N+1}},\xi ^{dp^{N+1}}}^{\left( h\right) }\left( \frac{%
a+idp^{N}}{dp^{N+1}}\right) \\
&&\hspace{-0.5in}=\left[ dp^{N}\right] _{q}^{n}\frac{\left[ 2\right] _{q}}{%
\left[ 2\right] _{q^{dp^{N}}}}\left( -1\right) ^{a}\xi ^{a}q^{ha}\left[ p\right] _{q^{dp^{N}}}^{n}\frac{\left[ 2\right] _{q^{dp^{N}}}%
}{\left[ 2\right] _{q^{dp^{N+1}}}}\sum_{i=0}^{p-1}\left( -1\right)
^{i}\left( \xi ^{dp^{N}}\right) ^{i}\left( q^{dp^{N}}\right)
^{hi}\\
&\times&E_{n,\left( q^{dp^{N}}\right) ^{p},\left( \xi
^{dp^{N}}\right)
^{p}}^{\left( h\right) }\left( \frac{\frac{a}{dp^{N}}+i}{p}\right) \\
&&\hspace{-0.5in}=\left[ dp^{N}\right] _{q}^{n}\frac{\left[ 2\right] _{q}}{%
\left[ 2\right] _{q^{dp^{N}}}}\left( -1\right) ^{a}\xi
^{a}q^{ha}E_{n,q^{dp^{N}},\xi ^{dp^{N}}}^{\left( h\right) }\left( \frac{a}{%
dp^{N}}\right)=\mu _{n,\xi ,q}^{\left( h\right) }\left( a+dp^{N}\mathbb{Z}%
_{p}\right) ,
\end{eqnarray*}%
\noindent where we have used (\ref{2.7}).

To present boundedness, we use equation (\ref{2.4}) to expand the polynomial
$E_{n,q^{dp^{N}},\xi ^{dp^{N}}}^{\left( h\right) }\left( \frac{a}{dp^{N}}%
\right) $, so that%
\begin{equation*}
\mu _{n,\xi ,q}^{\left( h\right) }\left(
a+dp^{N}\mathbb{Z}_{p}\right)=\frac{%
\left[ 2\right] _{q}}{\left( 1-q\right) ^{n}}\left( -1\right) ^{a}\xi
^{a}q^{ha}\sum_{j=0}^{n}\binom{n}{j}\left( -1\right) ^{j}q^{ja}\frac{1}{%
1+\xi ^{dp^{N}}q^{hdp^{N}+jdp^{N}}}.
\end{equation*}%
\noindent Now, since $d$ is an odd natural number and $p$ is an odd prime,
we have $\left| 1-\left( -\xi ^{dp^{N}}q^{hdp^{N}+jdp^{N}}\right) \right|
_{p}\geqslant 1$, so by induction on $j$, we obtain
\begin{equation*}
\left| \mu _{n,\xi ,q}^{\left( h\right) }\left( a+dp^{N}\mathbb{Z}%
_{p}\right) \right| _{p}\leqslant M
\end{equation*}
\noindent for a constant $M$. This is what we require, so the proof is
completed.
\end{proof}

Let $\chi $ be a Dirichlet character with conductor $d$. Then we can express
the generalized twisted $\left( h,q\right) $-Euler numbers associated with $%
\chi $ as an integral over $\mathbb{X}$, by using the measure $\mu _{n,\xi
,q}^{\left( h\right) }$.

\begin{lemma}
\label{lem1}For $n\in \mathbb{Z}$, $n\geqslant 0$, we have%
\begin{equation*}
\int\limits_{\mathbb{X}}\chi \left( t\right) d\mu _{n,\xi ,q}^{\left(
h\right) }\left( t\right) =E_{n,q,\xi ,\chi }^{\left( h\right) }.
\end{equation*}
\end{lemma}

\begin{proof}
From the definition of $p$-adic invariant integral, we have%
\begin{equation*}
\int\limits_{\mathbb{X}}\chi \left( t\right) d\mu _{n,\xi
,q}^{\left(
h\right) }\left( t\right)=\underset{N\rightarrow \infty }{\text{lim}}%
\sum_{c=0}^{dp^{N}-1}\chi \left( c\right) \left[ dp^{N}\right] _{q}^{n}\frac{%
\left[ 2\right] _{q}}{\left[ 2\right] _{q^{dp^{N}}}}\left(
-1\right) ^{c}\xi ^{c}q^{hc}E_{n,q^{dp^{N}},\xi ^{dp^{N}}}^{\left( h\right) }\left( \frac{c}{%
dp^{N}}\right) .
\end{equation*}%
\noindent Writing $c=a+dm$ with $a=0,1,\ldots ,d-1$ and $m=0,1,2,\ldots $,
we get%
\begin{eqnarray*}
\int\limits_{\mathbb{X}}\chi \left( t\right) d\mu _{n,\xi
,q}^{\left( h\right) }\left( t\right) &=&\left[ d\right] _{q}^{n}\frac{\left[ 2%
\right] _{q}}{\left[ 2\right] _{q^{d}}}\sum_{a=0}^{d-1}\chi \left( a\right)
\left( -1\right) ^{a}\xi ^{a}q^{ha} \\
&\times& \underset{N\rightarrow \infty }{\text{lim}}\left[
p^{N}\right] _{q^{d}}^{n}\frac{\left[ 2\right] _{q^{d}}}{\left[
2\right] _{\left( q^{d}\right)
^{p^{N}}}}\sum_{m=0}^{p^{N}-1}\left( -1\right) ^{m}\left( \xi
^{d}\right) ^{m}\left( q^{d}\right) ^{hm}E_{n,\left( q^{d}\right)
^{p^{N}},\left( \xi ^{d}\right) ^{p^{N}}}^{\left( h\right)
}\left( \frac{\frac{a}{d}+m}{p^{N}}\right) \\
&=&\left[ d\right] _{q}^{n}\frac{\left[ 2\right] _{q}}{\left[ 2\right]
_{q^{d}}}\sum_{a=0}^{d-1}\chi \left( a\right) \left( -1\right) ^{a}\xi
^{a}q^{ha}E_{n,q^{d},\xi ^{d}}^{\left( h\right) }\left( \frac{a}{d}\right) .
\end{eqnarray*}%
\noindent Assuming $\chi \left( 0\right) =0$ and by the fact that $\chi
\left( d\right) =0$, last expression equals $E_{n,\xi ,q,\chi }^{\left(
h\right) }$, and the proof is completed.
\end{proof}

Since it is impossible to have a non-zero translation-invariant measure on $%
\mathbb{X}$, $\mu _{n,\xi ,q}^{\left( h\right) }$ is not invariant under
translation, but satisfies the following:

\begin{lemma}
\label{lem2}For a compact-open subset $U$ of $\mathbb{X}$, we have%
\begin{equation*}
\mu _{n,\xi ,q}^{\left( h\right) }\left( pU\right) =\left[ p\right] _{q}^{n}%
\frac{\left[ 2\right] _{q}}{\left[ 2\right] _{q^{p}}}\mu _{n,\xi
^{p},q^{p}}^{\left( h\right) }\left( U\right) .
\end{equation*}
\end{lemma}

\begin{proof}
Let $U=a+dp^{N}\mathbb{Z}_{p}$ be the compact-open subset of $\mathbb{X}$.
Then%
\begin{eqnarray*}
\hspace{-0.3in}\mu _{n,\xi ,q}^{\left( h\right) }\left(
pU\right)&=&\mu
_{n,\xi ,q}^{\left( h\right) }\left( pa+dp^{N+1}\mathbb{Z}_{p}\right) \\
&=&\left[ dp^{N+1}\right] _{q}^{n}\frac{\left[ 2\right] _{q}}{\left[ 2\right]
_{q^{dp^{N+1}}}}\left( -1\right) ^{pa}\xi ^{pa}q^{hpa}E_{n,q^{dp^{N+1}},\xi
^{dp^{N+1}}}^{\left( h\right) }\left( \frac{pa}{dp^{N+1}}\right) \\
&=&\left[ p^{N}\right] _{q}^{n}\frac{\left[ 2\right] _{q}}{\left[ 2\right]
_{q^{p}}}\left[ dp^{N}\right] _{q^{p}}^{n}\frac{\left[ 2\right] _{q^{p}}}{%
\left[ 2\right] _{\left( q^{p}\right) ^{dp^{N}}}}\left( -1\right)
^{a}\left( \xi ^{p}\right) ^{a}\left( q^{p}\right)
^{ha}E_{n,\left( q^{p}\right)
^{dp^{N}},\left( \xi ^{p}\right) ^{dp^{N}}}^{\left( h\right) }\left( \frac{a%
}{dp^{N}}\right) \\
&=&\left[ p^{N}\right] _{q}^{n}\frac{\left[ 2\right] _{q}}{\left[ 2\right]
_{q^{p}}}\mu _{n,\xi ^{p},q^{p}}^{\left( h\right) }\left( a+dp^{N}\mathbb{Z}%
_{p}\right)=\left[ p\right] _{q}^{n}\frac{\left[ 2\right]
_{q}}{\left[ 2\right] _{q^{p}}}\mu _{n,\xi ^{p},q^{p}}^{\left(
h\right) }\left( U\right) ,
\end{eqnarray*}%
\noindent which is the desired result.
\end{proof}

Next, we give a relation between $\mu _{n,\xi ,q}^{\left( h\right) }$ and $%
\mu _{-q}$.

\begin{lemma}
\label{lem3}For any $n\in \mathbb{Z}$, $n\geqslant 0$, we have%
\begin{equation*}
d\mu _{n,\xi ,q}^{\left( h\right) }\left( t\right) =q^{\left( h-1\right)
t}\xi ^{t}\left[ t\right] _{q}^{n}d\mu _{-q}\left( t\right) .
\end{equation*}
\end{lemma}

\begin{proof}
From the definition of $\mu _{n,\xi ,q}^{\left( h\right) }$ and expansion of
twisted $\left( h,q\right) $-Euler polynomials, we have%
\begin{equation*}
\mu _{n,\xi ,q}^{\left( h\right) }\left(
a+dp^{N}\mathbb{Z}_{p}\right)=\frac{\left[ 2\right] _{q}}{\left(
1-q\right) ^{n}}\left( -1\right) ^{a}\xi
^{a}q^{ha}\sum_{j=0}^{n}\binom{n}{j}\left( -1\right) ^{j}q^{ja}\frac{1}{%
1+\xi ^{dp^{N}}q^{hdp^{N}+jdp^{N}}}.
\end{equation*}%
\noindent By the same method presented in \cite{Kim2002b}, we obtain%
\begin{eqnarray*}
\hspace{-0.3in}\underset{N\rightarrow \infty }{\text{lim}}\mu
_{n,\xi ,q}^{\left( h\right) }\left( a+dp^{N}\mathbb{Z}_{p}\right)
&=&\frac{1}{2}\frac{\left[ 2\right] _{q}}{\left( 1-q\right)
^{n}}\left( -1\right) ^{a}\xi
^{a}q^{ha}\sum_{j=0}^{n}\binom{n}{j}\left( -1\right) ^{j}q^{ja} \\
&=&\frac{1+q}{2}\xi ^{a}q^{\left( h-1\right) a}\left[ a\right]
_{q}^{n}\left( -1\right) ^{a}q^{a}=q^{\left( h-1\right) a}\xi ^{a}\left[ a\right] _{q}^{n}\underset{%
N\rightarrow \infty }{\text{lim}}\frac{\left( -1\right) ^{a}q^{a}}{\frac{%
1-\left( -q^{dp^{N}}\right) }{1-\left( -q\right) }} \\
&=&q^{\left( h-1\right) a}\xi ^{a}\left[ a\right] _{q}^{n}\underset{%
N\rightarrow \infty }{\text{lim}}\mu _{-q}\left( a+dp^{N}\mathbb{Z}%
_{p}\right) .
\end{eqnarray*}%
\noindent We thus have%
\begin{equation*}
d\mu _{n,\xi ,q}^{\left( h\right) }\left( t\right) =q^{\left( h-1\right)
t}\xi ^{t}\left[ t\right] _{q}^{n}d\mu _{-q}\left( t\right) ,
\end{equation*}%
\noindent the desired result.
\end{proof}

Let $\omega $ denote the Teichm\"{u}ller character mod $p$. For an arbitrary
character $\chi $ and $n\in \mathbb{Z}$, let $\chi _{n}=\chi \omega ^{-n}$
in the sense of product of characters. For $t\in \mathbb{X}^{\ast }=\mathbb{X%
}-p\mathbb{X}$, we set $\left\langle t\right\rangle _{q}=\left[ t\right]
_{q}/\omega \left( t\right) $. Since$\left| \left\langle t\right\rangle
_{q}-1\right| _{p}<p^{-1/\left( p-1\right) }$, $\left\langle t\right\rangle
_{q}^{s}$ is defined by exp$\left( s\text{log}_{p}\left\langle
t\right\rangle _{q}\right) $ for $\left| s\right| _{p}\leqslant 1$, where log%
$_{p}$ is the Iwasawa $p$-adic logarithm function (\cite{Iwasawa}). For $%
\left| 1-q\right| _{p}<p^{-1/\left( p-1\right) }$, we have $\left\langle
t\right\rangle _{q}^{p^{N}}\equiv 1\left( \text{mod}p^{N}\right) $.

We now define $p$-adic generalized twisted $\left( h,q\right) $-Euler-$l$%
-function.

\begin{definition}
\label{def1}For $s\in \mathbb{Z}_{p}$,%
\begin{equation*}
l_{p,q,\xi }^{\left( h\right) }\left( s,\chi \right) =\int\limits_{\mathbb{X}%
^{\ast }}\left\langle t\right\rangle _{q}^{-s}q^{\left( h-1\right) t}\xi
^{t}d\mu _{-q}\left( t\right) .
\end{equation*}
\end{definition}

\noindent The values of this function at non-positive integers are given by
the following theorem:

\begin{theorem}
\label{thm2}For any $n\in \mathbb{Z}$, $n\geqslant 0$,%
\begin{equation*}
l_{p,q,\xi }^{\left( h\right) }\left( -n,\chi \right) =E_{n,q,\xi ,\chi
_{n}}^{\left( h\right) }-\chi _{n}\left( p\right) \left[ p\right] _{q}^{n}%
\frac{\left[ 2\right] _{q}}{\left[ 2\right] _{q^{p}}}E_{n,q^{p},\xi
^{p},\chi _{n}}^{\left( h\right) }.
\end{equation*}
\end{theorem}

\begin{proof}
\begin{eqnarray*}
l_{p,q,\xi }^{\left( h\right) }\left( -n,\chi \right) &=&\int\limits_{%
\mathbb{X}^{\ast }}\left\langle t\right\rangle _{q}^{n}q^{\left( h-1\right)
t}\xi ^{t}d\mu _{-q}\left( t\right) =\int\limits_{\mathbb{X}^{\ast }}\chi
_{n}\left( t\right) \left[ t\right] _{q}^{n}q^{\left( h-1\right) t}\xi
^{t}d\mu _{-q}\left( t\right) \\
&=&\int\limits_{\mathbb{X}^{\ast }}\chi _{n}\left( t\right) d\mu
_{n,\xi ,q}^{\left( h\right) }\left(
t\right)=\int\limits_{\mathbb{X}}\chi _{n}\left( t\right) d\mu
_{n,\xi ,q}^{\left( h\right) }\left( t\right)
-\int\limits_{p\mathbb{X}}\chi _{n}\left( t\right)
d\mu _{n,\xi ,q}^{\left( h\right) }\left( t\right) \\
&=&E_{n,q,\xi ,\chi _{n}}^{\left( h\right) }-\chi _{n}\left( p\right) \left[
p\right] _{q}^{n}\frac{\left[ 2\right] _{q}}{\left[ 2\right] _{q^{p}}}%
E_{n,q^{p},\xi ^{p},\chi _{n}}^{\left( h\right) },
\end{eqnarray*}%
\noindent where Lemma \ref{lem1}, Lemma \ref{lem2} and Lemma \ref{lem3} are
used.
\end{proof}

This theorem will be mainly used in the next section, where certain
applications of $p$-adic generalized twisted $\left( h,q\right) $-Euler-$l$%
-function are given.

\section{\normalsize{Kummer Congruences for Generalized Twisted $\left( h,q\right) $%
-Euler Numbers}}

\setcounter{theorem}{0}

\hspace{0.15in}This section is devoted to an application of the
$p$-adic generalized twisted $\left( h,q\right)
$-Euler-$l$-function to an important number theoretic concept,
congruences systems. In particular, we derive Kummer-type
congruences for generalized twisted $\left( h,q\right) $-Euler
numbers by
using $p$-adic integral representation of $p$-adic generalized twisted $%
\left( h,q\right) $-Euler-$l$-function and Theorem \ref{thm2}.

In the sequel, we assume that $q\in \mathbb{C}_{p}$ with $\left| 1-q\right|
_{p}<1$. Then $q\equiv 1\left( \text{mod}\mathbb{Z}_{p}\right) $. For $t\in
\mathbb{X}^{\ast }$, we have $\left[ t\right] _{q}\equiv t\left( \text{mod}%
\mathbb{Z}_{p}\right) $, thus $\left\langle t\right\rangle _{q}\equiv
1\left( \text{mod}p\mathbb{Z}_{p}\right) $. For a positive integer $c$, the
forward difference operator $\Delta _{c}$ acts on a sequence $\left\{
a_{m}\right\} $ by $\Delta _{c}a_{m}=a_{m+c}-a_{m}$. The powers $\Delta
_{c}^{k}$ of $\Delta _{c}$ are defined by $\Delta _{c}^{0}=$identity and for
any positive integer $k$, $\Delta _{c}^{k}=\Delta _{c}\circ \Delta
_{c}^{k-1} $. Thus%
\begin{equation*}
\Delta _{c}^{k}a_{m}=\sum_{j=0}^{k}\binom{k}{j}\left( -1\right)
^{k-j}a_{m+jc}.
\end{equation*}%
\noindent For simplicity in the notation, we write%
\begin{equation*}
\varepsilon _{n,q,\xi ,\chi _{n}}^{\left( h\right) }=E_{n,q,\xi ,\chi
_{n}}^{\left( h\right) }-\chi _{n}\left( p\right) \left[ p\right] _{q}^{n}%
\frac{\left[ 2\right] _{q}}{\left[ 2\right] _{q^{p}}}E_{n,q^{p},\xi
^{p},\chi _{n}}^{\left( h\right) }.
\end{equation*}

\begin{theorem}
\label{thm3}For $n\in \mathbb{Z}$, $n\geqslant 0$ and $c\equiv 0\left( \text{%
mod}\left( p-1\right) \right) $, we have%
\begin{equation*}
\Delta _{c}^{k}\varepsilon _{n,q,\xi ,\chi _{n}}^{\left( h\right) }\equiv
0\left( \text{mod}p^{k}\mathbb{Z}_{p}\right) .
\end{equation*}
\end{theorem}

\begin{proof}
Since $\Delta _{c}^{k}$ is a linear operator, by Theorem \ref{thm2} we have%
\begin{eqnarray*}
\Delta _{c}^{k}\varepsilon _{n,q,\xi ,\chi _{n}}^{\left( h\right)
} &=&\Delta _{c}^{k}l_{p,q,\xi }^{\left( h\right) }\left( -n,\chi
\right)=\Delta _{c}^{k}\int\limits_{\mathbb{X}^{\ast
}}\left\langle t\right\rangle _{q}^{n}q^{\left( h-1\right) t}\xi
^{t}d\mu _{-q}\left(
t\right) \\
&=&\sum_{j=0}^{k}\binom{k}{j}\left( -1\right) ^{k-j}\int\limits_{\mathbb{X}%
^{\ast }}\left\langle t\right\rangle _{q}^{n+jc}q^{\left( h-1\right) t}\xi
^{t}d\mu _{-q}\left( t\right) \\
&=&\int\limits_{\mathbb{X}^{\ast }}\left\langle t\right\rangle
_{q}^{n}q^{\left( h-1\right) t}\xi ^{t}\left( \left\langle t\right\rangle
_{q}^{c}-1\right) ^{k}d\mu _{-q}\left( t\right) .
\end{eqnarray*}%
\noindent Now, $\left\langle t\right\rangle _{q}\equiv 1\left( \text{mod}p%
\mathbb{Z}_{p}\right) $, which implies that $\left\langle t\right\rangle
_{q}^{c}\equiv 1\left( \text{mod}p\mathbb{Z}_{p}\right) $ since $c\equiv
0\left( \text{mod}\left( p-1\right) \right) $, and thus%
\begin{equation*}
\left( \left\langle t\right\rangle _{q}^{c}-1\right) ^{k}\equiv 0\left(
\text{mod}p^{k}\mathbb{Z}_{p}\right) .
\end{equation*}%
\noindent Therefore%
\begin{equation*}
\Delta _{c}^{k}l_{p,q,\xi }^{\left( h\right) }\left( -n,\chi \right) \equiv
0\left( \text{mod}p^{k}\mathbb{Z}_{p}\right) ,
\end{equation*}%
\noindent from which the result follows.
\end{proof}

\begin{theorem}
\label{thm4}Let $n$ and $n^{\prime }$ be positive integers such that $%
n\equiv n^{\prime }\left( \text{mod}\left( p-1\right) \right) $. Then, we
have%
\begin{equation*}
\varepsilon _{n,q,\xi ,\chi _{n}}^{\left( h\right) }\equiv \varepsilon
_{n^{\prime },q,\xi ,\chi _{n^{\prime }}}^{\left( h\right) }\left( \text{mod}%
p\mathbb{Z}_{p}\right) .
\end{equation*}
\end{theorem}

\begin{proof}
Without loss of generality, let $n\geqslant n^{\prime }$. Then%
\begin{equation*}
l_{p,q,\xi }^{\left( h\right) }\left( -n,\chi \right) -l_{p,q,\xi }^{\left(
h\right) }\left( -n^{\prime },\chi \right) =\int\limits_{\mathbb{X}^{\ast
}}\left\langle t\right\rangle _{q}^{n}q^{\left( h-1\right) t}\xi ^{t}\left(
\left\langle t\right\rangle _{q}^{n-n^{\prime }}-1\right) d\mu _{-q}\left(
t\right) .
\end{equation*}%
\noindent Since $n-n^{\prime }\equiv 0\left( \text{mod}\left( p-1\right)
\right) $, we have $\left\langle t\right\rangle _{q}^{n-n^{\prime }}-1\equiv
0\left( \text{mod}p\mathbb{Z}_{p}\right) $, which entails the result.
\end{proof}

\bigskip

\textbf{Acknowledgment:} This work was supported by Akdeniz University
Scientific Research Projects Unit.

\end{document}